\documentclass[10pt]{article}

\usepackage{amssymb}
\usepackage{eucal}
\usepackage{amsmath}
\usepackage{amsthm}
\usepackage{amscd}
\usepackage{graphics}
\usepackage{graphicx}
\usepackage{amsfonts}
\usepackage{latexsym}
\usepackage[all]{xy}

\numberwithin{equation}{section}
\newtheorem{thm}{\bf Theorem}[section]
\newtheorem{lem}[thm]{\bf Lemma}
\newtheorem{prop}[thm]{\bf Proposition}

\newtheorem{corollary}[thm]{Corollary}
\theoremstyle{remark}
\newtheorem{rem}{\bf Remark}[section]

\begin{document}

\title{Instability conditions for circulatory and gyroscopic conservative systems}
\author{Petre Birtea\footnote{Corresponding author; West University of Timi\c soara, Department of Mathematics, Bd. Vasile P\^ arvan, No. 4, 300223 Timi\c soara, Romania; tel.: +40 726 424147, fax: +40 256 592316, e-mail: birtea@math.uvt.ro}, Ioan Ca\c{s}u, Dan Com\u anescu\\
{\small Department of Mathematics, West University of Timi\c soara,}\\
{\small Bd. Vasile P\^ arvan, No. 4, 300223 Timi\c soara, Romania}\\
{\small E-mail: birtea@math.uvt.ro; casu@math.uvt.ro; comanescu@math.uvt.ro}}
\date{}

\maketitle

\begin{abstract}
We give a method which generates sufficient conditions for instability of equilibria for circulatory and gyroscopic conservative systems. The method is based on the Gramians of a set of vectors whose coordinates are powers of the roots of the characteristic polynomial for the studied systems. New instability results are obtained for general circulatory and gyroscopic conservative systems. We also apply this method for studying the instability of motion for a charged particle in a stationary electromagnetic field.

\end{abstract}


\section{Introduction}

Many physical phenomenons are modeled by second order Euler-Lagrange equations
\begin{equation}\label{}
    \frac{d}{dt}\frac{\partial L}{\partial \dot{q}^i}-\frac{\partial L}{\partial q^i}=Q_i,
\end{equation}
where $L$ is the Lagrangian function and $Q_i$ are the components of generalized forces. The linearized equation at an equilibrium point has the form
\begin{equation}\label{linear}
    M\ddot{\bf q}+A_2\dot{\bf q}+A_3{\bf q}={\bf 0},
\end{equation}
where $M$ is a constant symmetric matrix and $A_2,A_3$ are constant matrices.
The above system appears, for example, in vibration theory and it is known as {\it lumped-mass system (without external forces)} (see \cite{inman}), where ${\bf q}$ is an $n$ vector of time-varying elements representing the displacements of the masses. Assuming that $M$ is invertible we multiply the equation \eqref{linear} by the matrix $M^{-1}$. Using the notations $D$ and $K$ for the symmetric parts, respectively $G$ and $C$ for the skew-symmetric parts of $M^{-1}A_2$ and $M^{-1}A_3$ we obtain the following normal system, see \cite{inman}, \cite{marsden}
\begin{equation}\label{lumped-normal}
    \ddot{\bf q}+(D+G)\dot{\bf q}+(K +C){\bf q}={\bf 0}.
\end{equation}
In literature, $D$ is called the {\it viscous damping matrix}, $G$ is the {\it gyroscopic matrix}, $K$ is the {\it stiffness matrix} and $C$ is the {\it circulatory matrix}.

According to \cite{inman} we have the following classification of the normal systems of the form \eqref{lumped-normal}:
\begin{itemize}
\item [(i)] {\it conservative systems} when $D=G=C=0$ and $K$ is positive definite;
\item [(ii)] {\it gyroscopic conservative systems} or {\it undamped gyroscopic systems} when $D=C=0$;
\item [(iii)] {\it damped non-gyroscopic systems} or {\it passive systems} when $G=C=0$ and $D,K$ are positive definite;
\item [(iv)] {\it circulatory systems} when $D=G=0$;
\item [(v)] {\it systems with constraint damping} when $G=0$.
\end{itemize}

In Section 2 we give a list of sufficient conditions for the existence of a complex root with non-zero imaginary part for a polynomial with real coefficients. We obtains these conditions using the Gramian of a set of vectors whose coordinates are powers of the roots for the given polynomial. We explicitly write three of these conditions, which turn out to be conditions involving only the power sums of the roots for the given polynomial. In the case of a characteristic polynomial of a matrix, the above mentioned conditions are written only in terms of norms and traces of the symmetric and skew-symmetric parts of the matrix.

In Section 3 we give sufficient conditions for the instability of circulatory systems. These conditions involve only the norms and traces of the stiffness matrix and the circulatory matrix. We recover a sufficient condition previously found by Bulatovic in \cite{bulatovic1}. We also give an example of a circulatory system, which shows that the three sufficient conditions presented in this section give different instability regions.

In Section 4 we give sufficient conditions for the instability of gyroscopic conservative systems involving norms and traces of the stiffness matrix and the gyroscopic matrix. We apply these conditions for the case of a charged particle acted by the Lorentz force. We also describe the regions of instability for this system.

\section{Sufficient conditions for complex roots with non-zero imaginary part}

Let $Q(\alpha)=\alpha^n+a_1\alpha^{n-1}+\dots +a_n$ be a $n$-degree polynomial with real coefficients. We denote by $\{\alpha_1,\dots,\alpha_n\}$ the set of roots of the polynomial $Q$. We search for sufficient conditions in order to have at least one complex root with non-zero imaginary part.

Consider the vectors:
\begin{equation*}
\left.\begin{array}{l}
v_0=(1,\dots,1);\\
v_1=(\alpha_1,\dots,\alpha_n);\\
\vdots \\
v_{n-1}=(\alpha_1^{n-1},\dots,\alpha_n^{n-1}).
\end{array}\right.
\end{equation*}

We have the well-known result: if all roots of the polynomial $Q$ are real, then for any $0\leq i_1\leq\dots\leq i_k\leq n-1$ the following inequality holds
$$\operatorname{Gram}(v_{i_1},\dots,v_{i_k}):=\left|\begin{array}{ccc}
<v_{i_1},v_{i_1}>&\dots&<v_{i_1},v_{i_k}>\\
\vdots&\ddots&\vdots\\
<v_{i_k},v_{i_1}>&\dots&<v_{i_k},v_{i_k}>\end{array}\right|\geq 0.$$
\begin{lem}
\label{lema}
If there exist $0\leq i_1<\dots < i_k\leq n-1$ such that
$$\operatorname{Gram}(v_{i_1},\dots,v_{i_k}) <0,$$
then the polynomial $Q$ has at least one complex root with non-zero imaginary part.
\end{lem}

The power sums of the roots are real numbers denoted by $s_k:=\sum\limits_{i=1}^n\alpha_i^k$, where $k\in {\mathbb N}$. We notice that
$$\operatorname{Gram}(v_{0},\dots,v_{{n-1}})=\left|\begin{array}{cccc}
n&s_1&\dots&s_{n-1}\\
s_1&s_2&\dots&s_n\\
\vdots&\vdots&\ddots&\vdots\\
s_{n-1}&s_n&\dots&s_{2n-2}\end{array}\right|.$$
Applying Lemma \ref{lema} for the simple particular cases $\operatorname{Gram}(v_0,v_1)$, and $\operatorname{Gram}(v_0,v_2)$ and $\operatorname{Gram}(v_1,v_2)$ we obtain the following sufficient conditions.

\begin{prop}
\label{prop1}
If any of the inequalities hold
\begin{itemize}
\item[(i)] $ns_2<s_1^2$;
\item[(ii)] $ns_4<s_2^2$;
\item[(iii)] $s_2s_4<s_3^2$,
\end{itemize}
then the polynomial $Q$ has at least one complex root with non-zero imaginary part.
\end{prop}

Considering Gramians of higher order, one can obtain similar sufficient conditions.
We notice that the inequalities (i) and (ii) hold when $s_2<0$, respectively $s_4<0$. Also, inequality (iii) holds, in particular, when $s_2<0$ and  $s_4>0$ or $s_2>0$ and  $s_4<0$.
If $P(\lambda)=Q(\lambda^2)$, then Proposition \ref{prop1} gives sufficient conditions for the existence of a complex root with strictly positive real part of the polynomial $P$.

All the power sums $s_k$ can be expressed recurrently, using Newton formulas, in terms of the coefficients of $Q$ only (see \cite{gantmaher}):
\begin{equation*}
\left.\begin{array}{l}
s_k+a_1s_{k-1}+\dots +a_{k-1}s_1+ka_k=0,~~~~~1\leq k\leq n\\
s_k+a_1s_{k-1}+\dots +a_{k-1}s_{k-n+1}+a_ns_{k-n}=0,~~~~~k> n.
\end{array}\right.
\end{equation*}

The above recurrences and Lemma \ref{lema} give sufficient conditions for the existence of complex roots with non-zero imaginary part exclusively in terms of the coefficients of the polynomial $Q$. More precisely, we have the following result.

\begin{prop}
\label{prop2}
The polynomial $Q$ has at least one complex root with non-zero imaginary part if:
\begin{itemize}
\item[(i)]
$$n\geq 2:~~~n(a_1^2-2a_2)<a_1^2;$$
\item[(ii)]
\begin{align*}
n=2:&~~~a_1^2(a_1^2-4a_2)<0;\\
n=3:&~~~a_1^4+6a_1a_3+a_2^2<4a_1^2a_2;\\
n\geq 4:&~~~n(a_1^4-4a_1^2a_2+4a_1a_3+2a_2^2-4a_4)<(a_1^2-2a_2)^2;
\end{align*}
\item[(iii)]
\begin{align*}
n=2:&~~~a_2^2(a_1^2-4a_2)<0;\\
n=3:&~~~a_1^2a_2^2+10a_1a_2a_3<2a_1^3a_3+4a_2^3+9a_3^2;\\
n\geq 4:&~~~a_1^2a_2^2+10a_1a_2a_3+8a_2a_4<2a_1^3a_3+4a_1^2a_4+4a_2^3+9a_3^2.
\end{align*}
\end{itemize}
\end{prop}

For the case when $Q$ is the characteristic polynomial of a matrix $M\in {\cal M}_{n\times n}({\mathbb R})$, i.e. $Q(\alpha)=\det (\alpha I_n-M)$, conditions (i), and (ii) and (iii) of Proposition \ref{prop1} can be interpreted only in terms of the symmetric and skew-symmetric parts of the matrix $M$.
Any matrix $M$ can be uniquely decomposed as a sum of a symmetric part and a skew-symmetric part, $M=M_s+M_a$, where $M_s=\frac{1}{2}(M+M^T)$  and
$M_a=\frac{1}{2}(M-M^T)$.

In what follows, we will use the well-known equalities: $\operatorname{Tr}(A+B)=\operatorname{Tr}(A)+\operatorname{Tr}(B)$, and $\operatorname{Tr}(AB)=\operatorname{Tr}(BA)$ for any $A,B\in {\cal M}_{n\times n}({\mathbb R})$ and $\operatorname{Tr}(AB)=0$ for any symmetric matrix $A$ and any skew-symmetric matrix $B$. Furthermore, we have $||A||^2=\operatorname{Tr}(AA^T)$ for any $A=[a_{ij}] \in {\cal M}_{n\times n}({\mathbb R})$, where $||A||^2:=\sum\limits_{i,j}a_{ij}^2$. Consequently, $||M_s||^2=\operatorname{Tr}(M_s^2)$ and $||M_a||^2=-\operatorname{Tr}(M_a^2)$.
By a direct computation we have
\begin{align*}
s_1&=\operatorname{Tr}(M)=\operatorname{Tr}(M_s);\\
s_2&=\operatorname{Tr}(M^2)=\operatorname{Tr}(M_s^2+M_a^2+M_sM_a+M_aM_s)\\
&=\operatorname{Tr}(M_s^2+M_a^2)=||M_s||^2-||M_a||^2;\\
s_3&=\operatorname{Tr}(M^3)=\operatorname{Tr}((M_s+M_a)^3)\\
&=\operatorname{Tr}(M_s^3+M_sM_a^2+M_s^2M_a+M_sM_aM_s
+M_aM_s^2+M_aM_sM_a+M_a^3+M_a^2M_s)\\
&=\operatorname{Tr}(M_s^3)+3\operatorname{Tr}(M_sM_a^2);\\
s_4&=\operatorname{Tr}(M^4)=\operatorname{Tr}((M_s+M_a)^4)\\
&=\operatorname{Tr}(M_s^4)+\operatorname{Tr}(M_a^4)+4\operatorname{Tr}(M_s^2M_a^2)+2\operatorname{Tr}((M_sM_a)^2)\\
&=||M_s^2||^2+||M_a^2||^2-4||M_sM_a||^2+2\operatorname{Tr}((M_sM_a)^2).
\end{align*}
The next theorem is a completion of a result given in \cite{bulatovic1} and it follows by substituting the above expressions for $s_1,s_2,s_3,s_4$ in Proposition \ref{prop1}.

\begin{thm}
\label{teorema}
If one of the inequalities hold
\begin{itemize}
\item[(i)] $n(||M_s||^2-||M_a||^2)<\operatorname{Tr}^2(M_s)$;
\item[(ii)] $n(||M_s^2||^2+||M_a^2||^2-4||M_sM_a||^2+2\operatorname{Tr}((M_sM_a)^2))<(||M_s||^2-||M_a||^2)^2$;
\item[(iii)] $(||M_s||^2-||M_a||^2)(||M_s^2||^2+||M_a^2||^2-4||M_sM_a||^2+2\operatorname{Tr}((M_sM_a)^2))<(\operatorname{Tr}(M_s^3)+3\operatorname{Tr}(M_sM_a^2))^2$,
\end{itemize}
then the matrix $M$ has at least one complex eigenvalue with non-zero imaginary part.
\end{thm}

\medskip

For $n=2$, if $\operatorname{Tr}(M)\not= 0$, then the inequalities (i) and (ii) from the above theorem are equivalent. As we will show later, this is not the case for $n\geq 3$.

\section{Instability for circulatory systems}
First we will give sufficient conditions of instability for a circulatory system subject to potential (conservative) and non-conservative positional (circulatory) forces
\begin{equation}
\label{circulatoriu}
\ddot {\bf q}+K{\bf q}+C{\bf q}=0,
\end{equation}
where $K$ is symmetric and $C$ is skew-symmetric. In \cite{bulatovic1} it has been showed that if the circulatory forces are bigger, in a certain sense, than the conservative forces, then the equilibrium point ${\bf q}={\bf 0},\dot {\bf q}={\bf 0}$ is unstable. The characteristic polynomial of the system \eqref{circulatoriu} is given by $P(\lambda)=\det(\lambda^2 I_n+K+C)$. This polynomial contains only even powers of $\lambda$. Denoting $\lambda^2=\alpha$ we obtain the characteristic polynomial $Q(\alpha)=\det(\alpha I_n+K+C)$ of the matrix $M=-(K+C)$. The existence of a complex root $\alpha$ with non-zero imaginary part implies the existence of a complex root $\lambda$ for the polynomial $P$ with strictly positive real part. Thus, instability of the system \eqref{circulatoriu} is implied by the existence of a complex root for the polynomial $Q(\alpha)$. Applying Theorem \ref{teorema} we obtain a different proof of a previous instability result presented in \cite{bulatovic1} and two new instability criteria.

\begin{thm}
\label{teorema2}
If one of the following inequalities hold
\begin{itemize}
\item[(i)] (\cite{bulatovic1}) $n(||K||^2-||C||^2)<\operatorname{Tr}^2(K)$;
\item[(ii)] $n(||K^2||^2+||C^2||^2-4||KC||^2+2\operatorname{Tr}((KC)^2))<(||K||^2-||C||^2)^2$;
\item[(iii)] $(||K||^2-||C||^2)(||K^2||^2+||C^2||^2-4||KC||^2+2\operatorname{Tr}((KC)^2))<(\operatorname{Tr}(K^3)+3\operatorname{Tr}(KC^2))^2$,
\end{itemize}
then the equilibrium point ${\bf q}={\bf 0},\dot {\bf q}={\bf 0}$ of the system \eqref{circulatoriu} is unstable.
\end{thm}

\begin{rem}
The inequality $ns_4<s_2^2$ is equivalent with $n\sum\limits_{i=1}^n (\alpha_i^2)^2<(\sum\limits_{i=1}^n\alpha_i^2)^2$. Noticing that $\alpha_i^2,~i=1,2,\dots,n$, are the roots for the characteristic polynomial $Q(\alpha)=\det (\alpha I_n-M^2)$ we can apply
the inequality (i) of Proposition \ref{prop1} to the polynomial $Q(\alpha)=\det (\alpha I_n-M^2)$. We obtain
$$||M^2_a||^2>||M^2_s||^2-\displaystyle\frac{1}{n}\operatorname{Tr}^2(M^2_s).$$
An elementary computation gives us that
$M^2_a=KC+CK$
and
$M^2_s=K^2+C^2.$
Consequently, the inequality (ii) of Theorem \ref{teorema2} is equivalent with
$$||CK+KC||^2>||K^2+C^2||^2-\displaystyle\frac{1}{n}(||K||^2-||C||^2)^2.~~~\bigtriangleup$$
\end{rem}

The inequalities $s_2<0$ or $s_4<0$ obviously imply the existence of at least one complex eigenvalue with non-zero imaginary part. As a consequence, we obtain the following sufficient conditions for instability:
\begin{corollary}
If one of the following inequalities hold
\begin{itemize}
\item[(i)] $||K||<||C||$;
\item[(ii)] $||K^2||^2+||C^2||^2+2\operatorname{Tr}((KC)^2)<4||KC||^2$,
\end{itemize}
then the equilibrium point ${\bf q}={\bf 0},\dot {\bf q}={\bf 0}$ of the system \eqref{circulatoriu} is unstable.
\end{corollary}

Other sufficient instability conditions for circulatory systems can also be found in \cite{bulatovic2}, \cite{gallina}.  

For $n=2$, if $\operatorname{Tr}(K)\not= 0$, then the inequalities (i) and (ii) from the above theorem are equivalent. For $n\geq 3$ the conditions (i), and (ii) and (iii) give different instability regions, as the following example shows. Let
$$K=\left[\begin{array}{ccc}
1&0&0\\
0&1&k\\
0&k&0\end{array}\right];~~~C=\left[\begin{array}{ccc}
0&c&0\\
-c&0&0\\
0&0&0\end{array}\right].$$
Condition (i) is equivalent with the inequality $3c^2-3k^2-1>0$, condition (ii) is equivalent with the inequality $-(c^2-k^2)^2+14c^2-2k^2-1>0$ and condition (iii) is equivalent with $4c^6-4k^6-12c^2k^2(c^2-k^2)+8c^4-3k^4+4c^2k^2+4c^2>0$. The Figure 1 below shows the three instability regions corresponding to conditions (i), (ii), respectively (iii).

\begin{figure}[!h]
\begin{center}
\includegraphics[scale=0.23,angle=0]{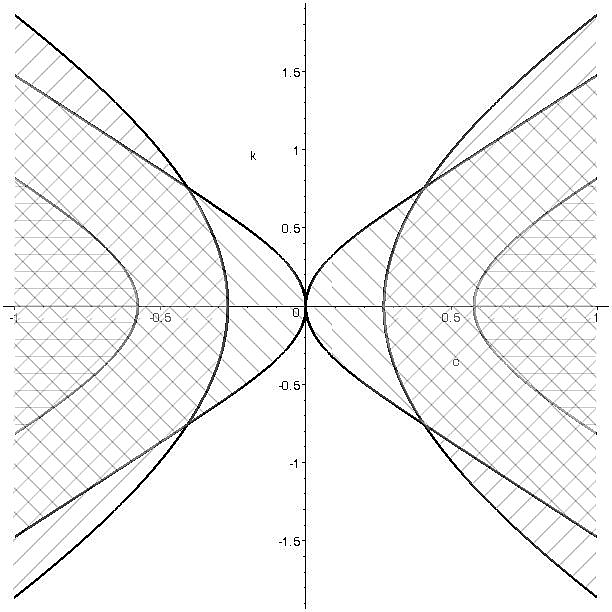}
\end{center}
\caption{\footnotesize
Instability regions using Theorem \ref{teorema2}: the horizontal grid designates the instability region obtained with condition (i); the right-inclined grid designates the instability region obtained with condition (ii); the left-inclined grid designates the instability region obtained with condition (iii).
}
\end{figure}

The above example shows that conditions (ii) and (iii) provide two supplementary instability regions that cannot be determined by using condition (i).

\section{Instability for gyroscopic conservative systems}

Next, we will give sufficient conditions of instability for a mechanical system subject to gyroscopic forces and potential forces
\begin{equation}
\label{gyro}
\ddot {\bf q}+G\dot {\bf q} +K{\bf q}={\bf 0},
\end{equation}
where $G$ is skew-symmetric and describes the gyroscopic forces and $K$ is symmetric and describes the potential forces. Various results concerning the stability problem for gyroscopic conservative systems can be found in \cite{huseyin}, \cite{yang}, \cite{afolabi}, \cite{bulatovic4}.

The matrix associated with the linear system \eqref{gyro} is given by
$$A=\left[\begin{array}{cc}
O_n&I_n\\
-K&-G
\end{array}\right]$$
and its characteristic equation is $P(\lambda):=\det(\lambda^2 I_n+\lambda G+K)=0$. Using the skew-symmetry of $G$ and the symmetry of $K$ we obtain that if $\lambda$ is a root of the characteristic equation, then $-\lambda$ is also a root of the characteristic equation. It follows that the above characteristic polynomial contains only powers of $\lambda^2$. After the substitution $\lambda^2=\alpha$ we obtained 
the reduced polynomial $Q(\alpha)$. 
We will apply Proposition \ref{prop1} and Proposition \ref{prop2} for the polynomial $Q(\alpha)$ in order to obtain sufficient conditions for the existence of complex roots with non-zero imaginary part. Thus, we get sufficient conditions for instability of the equilibrium point ${\bf q}={\bf 0},\dot {\bf q}={\bf 0}$ of the linear system \eqref{gyro}. 

We denote by $s_k^Q$ the power sums of the roots of the polynomial $Q(\alpha)$. We also denote by $s_k^P$ the power sums of the roots for the characteristic polynomial $P(\lambda)$. The following equalities hold $s_2^P=2s_1^Q$ and $s_4^P=2s_2^Q$. Consequently, the condition (i) of Proposition \ref{prop1} applied to the polynomial $Q(\alpha)$, $ns_2^Q<(s_1^Q)^2$ is equivalent with $2ns_4^P<(s_2^P)^2$, which is also the condition (ii) of the same proposition applied to the polynomial $P(\lambda)$. The inequality $ns_2^Q<(s_1^Q)^2$ implies that the polynomial $Q(\alpha)$ has at least one complex root with non-zero imaginary part and thus inequality $2ns_4^P<(s_2^P)^2$ implies that the polynomial $P(\lambda)$ has at least one root with strictly positive real part, which guarantees the instability of \eqref{gyro}.

As before, $s_2^P=\operatorname{Tr}(A^2)$ and $s_4^P=\operatorname{Tr}(A^4)$. 
We have the following computations

$$A_s=\frac{1}{2}(A+A^T)=\left[\begin{array}{cc}
O_n&\frac{1}{2}(I_n-K)\\
\frac{1}{2}(I_n-K)&O_n
\end{array}\right];~~~A_a=\frac{1}{2}(A-A^T)=\left[\begin{array}{cc}
O_n&\frac{1}{2}(I_n+K)\\
-\frac{1}{2}(I_n+K)&-G
\end{array}\right].$$
We denote by $B_1:=\frac{1}{2}(I_n-K),~~~B_2:=\frac{1}{2}(I_n+K)$.
Consequently, we have
\begin{align*}
A_s^2&=\left[\begin{array}{cc}
B_1^2&O_n\\
O_n&B_1^2
\end{array}\right];~~~
A_a^2=\left[\begin{array}{cc}
-B_2^2&-B_2G\\
GB_2&G^2-B_2^2
\end{array}\right];\\
A_s^4&=\left[\begin{array}{cc}
B_1^4&O_n\\
O_n&B_1^4
\end{array}\right];~~~A_a^4=\left[\begin{array}{cc}
-B_2G^2B_2+B_2^4& \ast\\
\ast&G^4-G^2B_2^2-B_2^2G^2-GB_2^2G+B_2^4
\end{array}\right];\\
A_s^2A_a^2&=\left[\begin{array}{cc}
-B_1^2B_2^2&\ast \\
\ast&B_1^2G^2-B_1^2B_2^2
\end{array}\right];~~~(A_sA_a)^2=\left[\begin{array}{cc}
(B_1B_2)^2& \ast \\
\ast&(B_1B_2)^2
\end{array}\right].
\end{align*}

Using the equalities $\operatorname{Tr}(G^2B_2^2)=\operatorname{Tr}(B_2^2G^2)=\operatorname{Tr}(GB_2^2G)=\operatorname{Tr}(B_2G^2B_2)$ we obtain
\begin{align*}
\operatorname{Tr}(A_s^4)&=\frac{1}{8}\left(n-4\operatorname{Tr}(K)+6\operatorname{Tr}(K^2)-4\operatorname{Tr}(K^3)+\operatorname{Tr}(K^4)\right);\\
\operatorname{Tr}(A_a^4)&=\operatorname{Tr}(G^4)+\frac{1}{8}\left(n+4\operatorname{Tr}(K)+6\operatorname{Tr}(K^2)+4\operatorname{Tr}(K^3)+\operatorname{Tr}(K^4)\right)\\
&-[\operatorname{Tr}(G^2)+2\operatorname{Tr}(G^2K)+\operatorname{Tr}(G^2K^2)];\\
4\operatorname{Tr}(A_s^2A_a^2)&=-\frac{n}{2}+\operatorname{Tr}(K^2)+\operatorname{Tr}(G^2)-2\operatorname{Tr}(KG^2)+\operatorname{Tr}(K^2G^2)-\frac{1}{2}\operatorname{Tr}(K^4);\\
2\operatorname{Tr}((A_sA_a)^2)&=\frac{1}{4}\left(n-2\operatorname{Tr}(K^2)+\operatorname{Tr}(K^4)\right). 
 \end{align*}

Using the equality $\operatorname{Tr}(G^2K)=\operatorname{Tr}(KG^2)$ we have
\begin{align*}
s_2^P&=\operatorname{Tr}(A_s^2)+\operatorname{Tr}(A_a^2)\\
 &=\frac{1}{2}\left(n-2\operatorname{Tr}(K)+\operatorname{Tr}(K^2)\right)+\operatorname{Tr}(G^2)-\frac{1}{2}\left(n+2\operatorname{Tr}(K)+\operatorname{Tr}(K^2)\right)\\
 &=-2\operatorname{Tr}(K)-||G||^2;\\
s_4^P&=\operatorname{Tr}(A_s^4)+\operatorname{Tr}(A_a^4)+4\operatorname{Tr}(A_s^2A_a^2)+2\operatorname{Tr}((A_sA_a)^2)\\
 &=2\operatorname{Tr}(K^2)+\operatorname{Tr}(G^4)-4\operatorname{Tr}(GKG)\\
 &=2\operatorname{Tr}(K^2)+\operatorname{Tr}(G^4)+4\operatorname{Tr}(G^TKG)\\
 &=2||K||^2+||G^2||^2+4\operatorname{Tr}(G^TKG).
 \end{align*}

Summarizing, we obtain the following instability result that was also achieved in \cite{bulatovic3} with a different argument.
\begin{thm}
\label{bul}
If the following inequality holds
$$2n(2||K||^2+||G^2||^2+4\operatorname{Tr}(G^TKG))<(2\operatorname{Tr}(K)+||G||^2)^2,$$
then the equilibrium point ${\bf q}={\bf 0},\dot {\bf q}={\bf 0}$ of the system \eqref{gyro} is unstable.
\end{thm}

For $n=3$ we will study the instability for the classical example of a charged particle in a stationary electromagnetic field. The particle is acted by the Lorentz force
${\bf F}_L=q{\bf E}(\mathbf{x})+q\dot{{\bf x}}\times {\bf B}(\mathbf{x})$, where $q$ is the electric charge of the particle, ${\bf E}(\mathbf{x})$ is the electric field, ${\bf B}(\mathbf{x})$ is the magnetic field and ${\bf x}$ is the position of the particle. From the stationarity assumption and Maxwell-Faraday equation we have that ${\bf E}(\mathbf{x})=\nabla \phi(\mathbf{x})$, where $\phi(\mathbf{x})$ is the electric potential. The equation of motion is $\ddot{{\bf x}}=\frac{q}{m}\nabla \phi({\bf x})+\frac{q}{m}\dot{{\bf x}}\times {\bf B}(\mathbf{x})$, where $m$ is the mass of the particle. At an equilibrium point $(\mathbf{x}_e,\mathbf{0})$ we have that $\nabla \phi(\mathbf{x}_e)=0$ and the linearized equation at this equilibrium point is
\begin{equation}
\label{linearmax}
\ddot{{\bf y}}=\frac{q}{m}\operatorname{Hess}\phi(\mathbf{x}_e){\bf y}+\frac{q}{m}\hat{\bf B}(\mathbf{x}_e)\dot{\bf y},
\end{equation}
where $\hat{\bf B}(\mathbf{x}_e)$ is the constant skew-symmetric $3\times 3$ matrix associated to the constant vector ${\bf B}(\mathbf{x}_e)$. By a convenient choice of coordinates the matrix $\displaystyle\frac{q}{m}\hat{\bf B}(\mathbf{x}_e)$ can be rendered in the canonical form
$$\frac{q}{m} \hat{\bf B}(\mathbf{x}_e)=\left[\begin{array}{ccc}
0&-c&0\\
c&0&0\\
0&0&0\end{array}\right].$$
We consider the potential function $\phi({\bf x})=\displaystyle\frac{\varepsilon}{2}x_1^2+\varepsilon x_2x_3$. Making the notations $G=-\frac{q}{m}\hat{\bf B}(\mathbf{x}_e)$ and $K=-\frac{q}{m}\operatorname{Hess}\phi(\mathbf{x}_e)$, the equation \eqref{linearmax} is in the form \eqref{gyro},
where
$$K=\left[\begin{array}{ccc}
k&0&0\\
0&0&k\\
0&k&0\end{array}\right];~~~G=\left[\begin{array}{ccc}
0&c&0\\
-c&0&0\\
0&0&0\end{array}\right],$$
with $k=-\frac{\varepsilon q}{m}$ and $c,k\not= 0$. The characteristic polynomial is $P(\lambda)=\lambda^6+(c^2+k)\lambda^4-k^2\lambda^2-k^3$. After substitution $\lambda^2=\alpha$, the reduced polynomial becomes $Q(\alpha)=\alpha^3+(c^2+k)\alpha^2-k^2\alpha-k^3$.
Condition (i) of Proposition \ref{prop2} applied for the polynomial $Q(\alpha)$ is equivalent with the condition from Theorem \ref{bul} and gives the inequality $-(c^2+k)^2-3k^2>0$. This inequality does not give an instability region.
Condition (ii) of Proposition \ref{prop2} applied for the polynomial $Q(\alpha)$ is equivalent with the inequality $c^6+4c^4k+10c^2k^2+6k^3<0$ and condition (iii) of Proposition \ref{prop2} applied for the polynomial $Q(\alpha)$ is equivalent with the inequalities $k<0$ and $2c^6+7c^4k+18c^2k^2+8k^3>0$.
The Figure 2 below shows the two instability regions corresponding to conditions (ii) and (iii) of Proposition \ref{prop2} applied to the polynomial $Q(\alpha)$.

\begin{figure}[!h]
\begin{center}
\includegraphics[scale=0.25,angle=0]{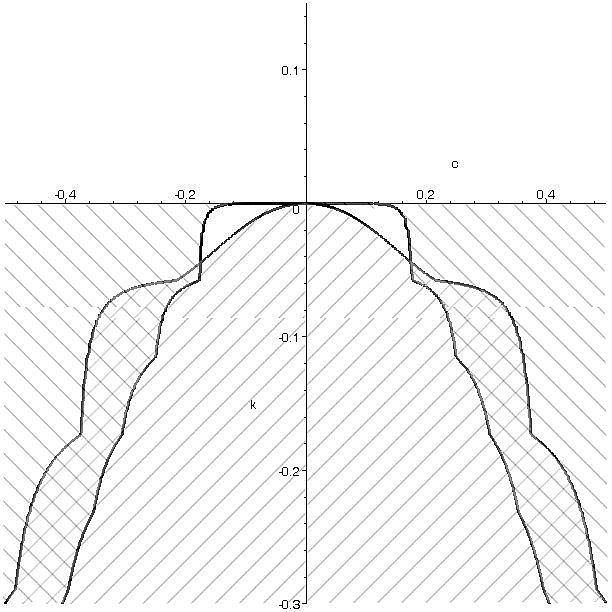}
\end{center}
\caption{\footnotesize
Instability regions using Proposition \ref{prop2}: the right-inclined grid designates the instability region obtained with condition (ii); the left-inclined grid designates the instability region obtained with condition (iii).
}
\end{figure}

For $n = 2$ the characteristic polynomial can be expressed in terms of traces
and determinants of the matrices which appear in the normal form \eqref{lumped-normal}. The
study of the instability region (when the matrices $D$ and $C$ are also present) has been completely
solved in \cite{kirillov}, \cite{kirillov2010}, using a result from
\cite{bottema}.  

\section{Conclusions}

In this paper we have used the Gramian technique for obtaining sufficient conditions for the existence of at least one complex root with non-zero imaginary part for a polynomial. Applying these results we have obtained sufficient instability conditions for circulatory systems and gyroscopic conservative systems. In this way, we recover the main theorems in \cite{bulatovic1} and \cite{bulatovic3}. Using other Gramians of second or higher order, one can obtain similar sufficient instability conditions. Gramians have also been used in studying the geometry of a higher order dissipation \cite{bico}.

\medskip

\noindent {\bf Acknowledgements.} This work was supported by a grant of the Romanian National Authority for Scientific Research, CNCS – UEFISCDI, project number PN-II-RU-TE-2011-3-0006.

\end{document}